\documentclass[a4paper,12pt]{article}
\usepackage{amssymb}
\usepackage{amsmath}
\usepackage{amsthm}
\usepackage{graphicx}
\numberwithin{equation}{section}
\newcommand{\ep}{\varepsilon}

\newcommand{\va}{\varphi}
\newcommand{\ppp}{\partial}

\newcommand{\weight}{e^{2s\va}}
\newcommand{\sumj}{\sum_{j=1}^n}

\newcommand{\R}{\mathbb{R}}

\newcommand{\N}{\mathbb{N}}

\newcommand{\www}{\widetilde}

\newcommand{\ooo}{\overline}
\newcommand{\OOO}{\Omega}

\newcommand{\sumd}{\sum_{k=1}^d}
\allowdisplaybreaks

\title
{A Carleman estimate and an energy method for a first-order symmetric 
hyperbolic system}

\author{
$^1\;$Giuseppe Floridia
$^2\;$Hiroshi Takase
$^{3.4,5,6}\;$Masahiro Yamamoto }

\begin{document}
\maketitle
\thanks{
$^1\;$ 
Mediterranean University of Reggio Calabria, 
Department PAU, \\
Via dell'Universit\`a 25  
89124 Reggio Calabria, Italy \\
INdAM Unit, University of Catania, Italy\\
e-mail:{\tt floridia.giuseppe@icloud.com}
\\
$^2\;$ Graduate School of Mathematical Sciences, The University
of Tokyo, Komaba, Meguro, Tokyo 153-8914, Japan \\
e-mail:{\tt htakase@ms.u-tokyo.ac.jp}
\\
$^3\;$Graduate School of Mathematical Sciences, The University
of Tokyo, Komaba, Meguro, Tokyo 153-8914, Japan \\
$^4\;$Honorary Member of Academy of Romanian Scientists, 
Ilfov, nr. 3, Bucuresti, Romania \\
$^5$ Correspondence member of Accademia Peloritana dei Pericolanti\\
$^6\;$Peoples' Friendship University of Russia 
(RUDN University) 6 Miklukho-Maklaya St, Moscow, 117198, Russian Federation
e-mail: {\tt myama@ms.u-tokyo.ac.jp}
}

\begin{abstract}
For a symmetric hyperbolic system of the first order, we prove 
a Carleman estimate under some positivity condition concerning  
the coefficient matrices.  Next, applying the Carleman estimate,
we prove an observability $L^2$-estimate for initial values by boundary data.
\end{abstract}

\section{Introduction}

Let $\OOO \subset \R^d$ be a bounded domain with smooth boundary 
$\ppp\OOO$, and let $N \in \N$ and 
$$
Q := \OOO \times (0,T).
$$
For any $k\in \{0,1,2,..., d\}$, let
$H^k = H^k(x,t)$ be $N\times N$ matrices whose elements belong to 
$C^1(\ooo{Q})$, and $P$ be an $N\times N$ matrix with elements in 
$L^{\infty}(Q)$.
We write $x= (x_1,..., x_d) \in \R^d$, and 
$$
\ppp_0 := \frac{\ppp}{\ppp t}, \quad
\ppp_k:= \frac{\ppp}{\ppp x_k}, \quad k=1,..., d.
$$
Here and henceforth $\cdot^T$ means the transpose of vectors and 
matrices under consideration. 

Throughout this article, we 
assume that $H^k(x,t)$, $k=0,1,..., d$, are symmetric matrices,
that is, 
$$
(H^k)^T = H^k, \quad k=0,1,..., d.
$$ 
Let an $\R^N$-valued function $F(x,t)$ be in $L^2(Q)$.
 
We consider a symmetric hyperbolic system of the first order:
$$
Lu(x,t) := H^0(x,t)\ppp_tu + \sum_{k=1}^d H^k(x,t)\ppp_ku 
+ P(x,t)u = F(x,t), \quad x\in \OOO, \, 0<t<T,
                                   \eqno{(1.1)}
$$
where $u = (u_1, ..., u_N)^T$.

Whenever there is no fear of confusion, we understand 
that $(U\cdot V)$ denotes a scalar product in $\R^N$ or
$L^2(\OOO)$ for $U=(U_1, ..., U_N)^T$ and $V = (V_1,..., V_N)^T$:
$$
(U\cdot V) := \sum_{j=1}^N U_jV_j \quad\mbox{or}\quad
(U\cdot V):= \int_{\OOO} \sum_{j=1}^N U_j(x)V_j(x) dx.
$$
By $\Vert \cdot\Vert_X$ we denote the norm in a normed space $X$
for specifying $X$, and let $\Vert \cdot\Vert$ mean 
$\Vert \cdot\Vert_{L^2(\OOO)}$ if we do not 
specify extra, and $\vert U\vert_{\R^N}:= \left( \sum_{j=1}^N U_j^2\right)
^{\frac{1}{2}}$.
\\

The system (1.1) can describe several important equations in mathematical
physics such as Maxwell's equations and the elasticity equations.

The main purpose of this articile is to establish a Carleman estimate for the 
symmetric hyperbolic system (1.1), and apply it to an energy estimate called 
an observability inequality.

A Carleman estimate is a weighted $L^2$-estimate for solutions to (1.1), and
originates from Carleman \cite{Ca}, and provides an effective method to 
the unique continuation problem and inverse problems for partial differential 
equation.  As for the applications to the unique continuation, we can
refer to H\"ormander \cite{Ho}, Isakov \cite{Is} for example.
On the other hand, as pioneerng work for an application of a Carleman 
estimate to inverse problems, see Bukhgeim and Klibanov \cite{BK}.
See also Beilina and Klibanov \cite{BeKl}, Bellassoued and Yamamoto 
\cite{BY}, Imanuvilov and Yamamoto \cite{IY1}, \cite{IY2},
Klibanov \cite{Kl}, Klibanov and Timonov \cite{KT}, Yamamoto \cite{Ya}.
Here we do not intend any comprehensive lists of works, and the readers can 
consult also the references therein.

We emphasize that Carleman estimates are essential ingredients 
for solving the unique continuation 
and mathematically analyzing inverse problems.

For $N=1$, as related works by Carleman estimates, we refer to 
Gaitan and Ouzzane \cite{GO}, Cannarsa, Floridia, G\"olgeleyen, and 
Yamamoto \cite{CFGY}, Cannarsa, Floridia and Yamamoto \cite{CFY},
Floridia and Takase \cite{FT1}, \cite{FT2}, G\"olgeleyen and 
Yamamoto \cite{GY}, and see also 
Klibanov and Pamyatnykh \cite{KP1}, \cite{KP2}, Klibanov and 
Yamamoto \cite{KY}, Machida and Yamamoto \cite{MY} 
as for related problems for 
the radiative transport equation called the Boltzmann equation.
On the other hand, for a symmteric hyperbolic system of the first order, the 
research for the direct problem is completed (e.g., Mizohata \cite{Mi},
Petrovsky \cite{Pe}).  However, to the best knowledge of the authors, 
Carleman estimates are not available for (1.1) and so works for inverse problems by 
Carleman esimates have not been published.
\\

Now for the statements of our main results, we introduce some notations and
conditions.
Let $\nu = (\nu_1, ..., \nu_d)^T$ be the unit outward normal vector to 
$\ppp\OOO$ at $x$.  For $0 \le t \le T$, we set 
$$
\ppp\OOO_+(t) := \left\{ x\in \ppp\OOO;\, 
\left(\left(\sumd H^k(x,t)\nu_k\right)v \,\cdot \, v\right)
> 0 \quad \mbox{for all $v\in \R^N$}\right\}
                                                         \eqno{(1.2)}
$$
and
$$
\ppp\OOO_-(t) := \left\{ x\in \ppp\OOO;\, 
\left(\left(\sumd H^k(x,t)\nu_k\right)v \,\cdot \, v\right)
\le 0 \quad \mbox{for all $v\in \R^N$}\right\}.
                                                         \eqno{(1.3)}
$$
For $t$-independent $H^k$, we note that 
$\ppp\OOO_+(t)$ and $\ppp\OOO_-(t)$ are independent of $t$.
We remark that $\ppp\OOO \supsetneqq \ppp\OOO_+(t) \cup
\ppp\OOO_-(t)$, in general.
\\

Let $\va = \va(x,t) \in C^1(\ooo{Q})$.
Now we can state our first main result.
\\
{\bf Theorem 1 (Carleman estimate for a symmetric hyperbolic system of 
the first order)}.\\
{\it
We assume that there exists a constant $\delta>0$ such that 
$$
\left( \left( \sum_{k=0}^d (\ppp_k\va)H^k(x,t)\right)v\,\cdot\,v\right)
\ge \delta\vert v\vert_{\R^N}^2
\quad \mbox{for all $(x,t)\in Q$ and $v\in \R^N$}.
                                                      \eqno{(1.4)}
$$
Then there exist constants $C>0$ and $s_0>0$ such that 
$$
s\int_{\OOO} (H^0(x,0)u(x,0)\,\cdot\, u(x,0)) e^{2s\va(x,0)} dx
+ s^2\int_Q \vert u\vert^2e^{2s\va} dxdt               \eqno{(1.5)}
$$
$$
+ s\int^T_0\int_{\ppp\OOO_-(t)} \left\vert
\left(\left(\sumd H^k(x,t)\nu_k\right)u\,\cdot \,u\right)\right\vert
\weight dSdt
$$
$$
\le C\int_Q \vert F\vert^2 \weight dxdt 
+ Cs\int^T_0\int_{\ppp\OOO \setminus \ppp\OOO_-(t)} 
\vert u\vert^2 \weight dSdt
$$
$$
+ s\int_{\OOO} (H^0(x,T)u(x,T)\,\cdot\, u(x,T)) e^{2s\va(x,T)} dx
$$
for all $s \ge s_0$ and $u\in H^1(Q)$ satisfying (1.1) in $Q$.
}
\\

Theorem 1 is a direct generalization to the case of a system, where 
the case $N=1$ is considered 
by \cite{CFGY} and \cite{GY}, and the
weight function $\va$ is different from \cite{GO} and is of the same type
as \cite{CFGY}, \cite{GY}, \cite{MY}.
In general, Carleman estimates do not estimate $u(\cdot,0)$, but 
the left-hand side of (1.5) can include the values $u(x,0)$.

The general theory for the Carleman estimate was designed mostly 
for functions with compact supports (e.g., \cite{Ho}, \cite{Is}), 
while our Carleman 
estimate does not assume that $u$ has compact supports, which allows us to 
apply a simplified argument by Huang, Imanuvilov and Yamamoto \cite{HIY}.
If we apply Carleman estimates with compact supports, then we need a cut-off 
function, which is a quite conventional way.
However, the cut-off procedure makes the total arguments more complicated.
The application of 
a general theory for constructing Carleman estimates for functions without
compact supports, is also very complicated and the direct way is more relevant 
for the derivation of a Carleman estimate for (1.1).

In the case where $H^k$, $k=0, 1, ..., d$ are diagonal matrices,
system (1.1) is not coupled with the first-order terms but only with the 
zeroth order terms.  Then under adequate conditions, thanks to the large 
parameter $s>0$, the Carleman estimate can be directly derived from the case 
for $N=1$, so that our main interest in Theorem 1 is for the case of
non-diagonal matrices $H^k$ for $0\le k \le d$.
\\

Next we apply Theorem 1 to the observability inequality of 
estimating initial value by boundary data.
First in Theorem 1, we choose the weight function $\va(x,t)$ by
$$
\va(x,t) = \eta(x) - \beta t,
$$
where $\beta > 0$ and $\eta \in C^1(\ooo{\OOO})$ are chosen later.
We assume that there exist constants $\delta_0>0$, $\delta_1>0$ and
$M>0$ such that 
$$
\left( \sumd (\ppp_k\eta)(x)H^k(x,t)v\, \cdot\, v\right)
\ge \delta_0\vert v\vert_{\R^N}^2 \quad
\mbox{for all $(x,t) \in Q$ and $v\in \R^N$}.     \eqno{(1.6)}
$$
and
$$
\delta_1 \vert v\vert_{\R^N}^2 \le 
(H^0(x,t)v\, \cdot\, v) \le M\vert v\vert^2_{\R^N} \quad
\mbox{for all $(x,t) \in Q$ and $v\in \R^N$}.     \eqno{(1.7)}
$$
\\
{\bf Theorem 2 (Observability inequality).}
\\
{\it
We further assume (1.6), (1.7) and 
$$
T > \frac{M}{\delta_0}
\left( \max_{x\in \ooo{\OOO}} \eta(x) - \min_{x\in \ooo{\OOO}} \eta(x)
\right).                               \eqno{(1.8)}
$$
Then there exists a constant $C>0$ such that 
$$
\Vert u(\cdot,0)\Vert \le C\Vert u\Vert_{L^2(\ppp\OOO\times (0,T))}
$$
for all $u \in H^1(Q)$ satisfying 
$H^0(x,t)\ppp_tu + \sumd H^k(x,t)\ppp_ku + Pu = 0$ in $Q$.
}
\\

This article is composed of three sections.  In Section 2 we prove 
Theorem 1 and Section 3 is devoted to the proof of Theorem 2.
\section{Proof of Theorem 1}

The proof is based on integration by parts, which is applicable also for
direct proofs of Carleman estimates for 
parabolic and hyperbolic equations.  The main steps are described as 
follows.
\\
{\bf Step I.}
We transform system (1.1) in $u$ to a system in terms of 
$w:= ue^{s\va}$.
\\
{\bf Step II.}
Calculating the squared $L^2$-norms of the transformed system, 
we make a lower estimate of the cross terms composed of products 
of $w$ and $\ppp_kw$ for $k=0,1,..., d$.
\\
{\bf Step III.}
The lower estimate in Step II cannot dierctly produce 
the Carleman estimate, which is the same for parabolic and hypebolic
equations.  Therefore, we take scalar products of the transformed system by 
$B(x,t)w$ with suitably chosen matrix function $B(x,t)$ and we add
the resulting equality to the lower estimate obtained in Step II, which 
leads us to the conclusion.
\\

Before starting the proof, we remark that it is sufficient to prove 
Theorem 1 for $P\equiv 0$.  Indeed let Theorem 1 be proved in the case of
$P\equiv 0$.
Then
\begin{align*}
& s\int_{\OOO} (H^0(x,0)u(x,0)\,\cdot\, u(x,0)) e^{2s\va(x,0)} dx
+ s^2\int_Q \vert u\vert^2e^{2s\va} dxdt\\
+ & s\int^T_0\int_{\ppp\OOO_-(t)} \left\vert
\left(\left(\sumd H^k(x,t)\nu_k\right)u\,\cdot \,u\right)\right\vert
\weight dSdt\\
\le & C\int_Q \vert F-Pu\vert^2 \weight dxdt 
+ \int^T_0\int_{\ppp\OOO \setminus \ppp\OOO_-(t)} 
\vert u\vert^2 \weight dSdt\\
+ & s\int_{\OOO} (H^0(x,T)u(x,T)\,\cdot\, u(x,T)) e^{2s\va(x,T)} dx.
\end{align*}
Since 
$$
\vert F-Pu\vert^2 \le 2\vert F\vert^2 + 2\vert Pu\vert^2
\le 2\vert F\vert^2 + 2\Vert P\Vert_{L^{\infty}(Q)}^2\vert u\vert^2
$$
in $Q$, we can estimate $\int_Q \vert F-Pu\vert^2 \weight dxdt$ by 
$$
C_0\int_Q \vert F\vert^2 \weight dxdt + C_0\int_Q \vert u\vert^2
\weight dxdt
$$
with some constant $C_0 > 0$.  Choosing $s>0$ large, we can absorb the term  
$C_0\int_Q \vert u\vert^2 \weight dxdt$ into 
$s^2\int_Q \vert u\vert^2 \weight dxdt$ on the left-hand side.
Therefore, (1.5) with $P\in L^{\infty}(Q)$ is derived from (1.5) in the case
of $P\equiv 0$.
\\

Moreover we show
$$
(R(x,t)\ppp_kw\, \cdot\, w)_{\R^N}
= \frac{1}{2}\ppp_k(Rw\, \cdot\, w)_{\R^N}
- \frac{1}{2}(\ppp_kRw\cdot w)_{\R^N}, \quad k=0,1,..., d \quad \mbox{in $Q$},
                                                       \eqno{(2.1)}
$$
where $R \in C^1(\ooo{Q})$ is an $N\times N$ symmetric matrix and 
$w = (w_1,..., w_N)^T \in H^1(Q)$.

Indeed, by the symmetry of $R$, we have
$(R\ppp_kw\cdot w)_{\R^N} = (\ppp_kw \cdot Rw)_{\R^N}$, and so 
\begin{align*}
& \ppp_k(Rw\cdot w)_{\R^N} 
= (R\ppp_kw\cdot w)_{\R^N} + (Rw\cdot \ppp_kw)_{\R^N} 
+ ((\ppp_kR)w\cdot w)_{\R^N} \\
=& 2(R\ppp_kw\cdot w)_{\R^N} + ((\ppp_kR)w\cdot w)_{\R^N} \quad \mbox{in $Q$},
\end{align*}
which verifies (2.1).
\\
\vspace{0.2cm}
{\bf Step I.}

Not only for functions in $L^2(\OOO)$, but also 
for $u = (u_1, ..., u_N)^T, \, v=(v_1, ..., v_N)^T \in L^2(Q)$, we use the 
notation  
$$
(u\cdot v) := \int_Q \sum_{j=1}^N u_j(x,t)v_j(x,t) dxdt.
$$
We set 
$$
Lu:= H^0(x,t)\ppp_tu + \sumd H^k(x,t)\ppp_ku = F \quad \mbox{in $Q$}. 
                        \eqno{(2.2)}
$$
We set 
$$
w := e^{s\va}u, \quad Pw:= e^{s\va}L(e^{-s\va}w).
$$
Then 
$$
Pw = H^0\ppp_tw + \sumd H^k\ppp_kw
- s\left( \sumd (\ppp_k\va)H^k + (\ppp_t\va)H^0\right)w
$$
in $Q$.  For short description, we set
$$
A:= \sum_{k=0}^d (\ppp_k\va)H^k \quad \mbox{in $Q$},
$$
where $\ppp_0\va = \ppp_t\va$.  Then we can write (1.1) in terms of
the transformed system in terms of $w$:
$$
Pw = H^0\ppp_tw + \sumj H^k\ppp_kw - sAw = Fe^{s\va} \quad 
\mbox{in $Q$}.                              \eqno{(2.3)}
$$

{\bf Step II.}

By the definition, we have
$$
\Vert Pw\Vert^2 = \int_Q \vert Lw\vert^2 \weight dxdt 
= \int_Q \vert F\vert^2 \weight dxdt.
$$
Therefore, we have to make a lower estimate of $\Vert Pw\Vert^2$ and 
extract the terms on the left-hand side of (1.5).   
Direct calculations yield
\begin{align*}
& \Vert Pw\Vert^2 
= 2\left(\left(H^0\ppp_tw + \sumd H^k\ppp_kw\right)\, \cdot \,  -sAw\right)
+ s^2\Vert Aw\Vert^2 + \left\Vert H^0\ppp_tw + \sumd H^k\ppp_kw\right\Vert^2\\
\ge& -2s\int_Q \left(\sum_{k=0}^d H^k\ppp_kw \cdot Aw\right) dxdt
+ s^2\Vert Aw\Vert^2.
\end{align*}
Hence,
$$
-2s\int_Q \left(\sum_{k=0}^d H^k\ppp_kw \cdot Aw\right) dxdt
+ s^2\Vert Aw\Vert^2 \le \Vert Fe^{s\va}\Vert^2.           \eqno{(2.4)}
$$
The last term $s^2\Vert Aw\Vert^2$ can estimate $s^2\Vert w\Vert^2$ 
suitably, but we cannot directly estimate the cross term
$$
-2s\int_Q \left(\sum_{k=0}^d H^k\ppp_kw \cdot Aw\right) dxdt.
$$ 
The difficulty comes from the non-symmetry in the form 
$\left(\sum_{k=0}^d H^k\ppp_kw \cdot Aw\right)$, so that we cannot 
treat the derivatives $\ppp_kw$.
We remark that $H^k$ and $A$ are symmetric but 
$AH^k$ and $H^kA$ are not necessarily symmetric.
Thus we need some argument in Step III.

{\bf Step III.}

For compensating for the above cross term, we take the scalar product 
of (2.3) with $B(x,t)w$, where we choose $B$ later:
$$
\left( \sum_{k=0}^d H^k\ppp_kw \, \cdot\, Bw\right)
- s(Aw\cdot Bw) = (Fe^{s\va} \cdot Bw).
$$
For adjusting the power of $s$, we multiply by $2s$:
$$
2s\left( \sum_{k=0}^d H^k\ppp_kw \, \cdot\, Bw\right)
- 2s^2(Aw\cdot Bw) = 2s(Fe^{s\va} \cdot Bw).              \eqno{(2.5)}
$$
Adding (2.4) and (2.5), we obtain
\begin{align*}
& 2s\left( \sum_{k=0}^d H^k\ppp_kw \, \cdot\, (B-A)w\right)
+ s^2(\Vert Aw\Vert^2 - 2(Aw\cdot Bw))\\
\le & \Vert Fe^{s\va}\Vert^2 + 2s(Fe^{s\va} \cdot Bw).
\end{align*}
Applying the Cauchy-Schwarz inequality to the right-hand side and 
using
$$
2s(Fe^{s\va} \cdot Bw) 
= 2\left( \frac{Fe^{s\va}}{\ep} \, \cdot \, s\ep Bw\right)
\le \frac{1}{\ep^2}\Vert Fe^{s\va}\Vert^2 + s^2\ep^2\Vert Bw\Vert^2,
$$
we obtain
$$
2s\left( \sum_{k=0}^d H^k\ppp_kw \, \cdot\, (B-A)w\right)
+ s^2(\Vert Aw\Vert^2 - 2(Aw\cdot Bw))
$$
$$
\le \left( 1 + \frac{1}{\ep^2}\right)\Vert Fe^{s\va}\Vert^2 
+ s^2\ep^2\Vert Bw\Vert^2.                \eqno{(2.6)}
$$
In order to apply (2.1) to the first term on the left-hand side of (2.6)
for estimating the cross term, we choose $B$ suitably, for example, such that 
$$
B-A = -E_N,
$$
where $E_N$ is the $N\times N$ identity matrix.
Then  
$$
-2s\left( \sum_{k=0}^d H^k\ppp_kw \, \cdot\, w\right)
+ s^2(2(Aw\, \cdot \,w) - \Vert Aw\Vert^2 )
$$
$$
\le \left( 1 + \frac{1}{\ep^2}\right)\Vert Fe^{s\va}\Vert^2 
+ s^2\ep^2\Vert (A-E_N)w\Vert^2.                \eqno{(2.7)}
$$

With a parameter $\mu>0$, we set 
$$
\www\va:= \mu\va, \quad 
\www A(x,t) := \sum_{k=0}^d (\ppp_k\www\va)H^k = \mu A    \quad \mbox{in $Q$}.
$$

Henceforth $C>0$ denotes generic constants which are independent of 
$s>0$ and $\mu>0$.  
Setting $\www w:= ue^{s\www\va}$, we repeat the calculations to 
obtain
$$
-2s\left( \sum_{k=0}^d H^k\ppp_k\www w \, \cdot\, \www w\right)
+ s^2(2(\www A\www w\, \cdot \,\www w) - \Vert \www A\www w\Vert^2 )
$$
$$
\le \left( 1 + \frac{1}{\ep^2}\right)\Vert Fe^{s\www\va}\Vert^2 
+ s^2\ep^2\Vert (\www A-E_N)\www w\Vert^2.                \eqno{(2.8)}
$$
We can directly verify 
$$
\Vert \www A \www w\Vert^2 = \mu^2\Vert A \www w\Vert^2 
\le C\mu^2\Vert \www w\Vert^2
$$
and
$$
\Vert (\www A-E_N)\www w\Vert^2 \le 2\Vert \www A\www w\Vert^2 
+ 2\Vert E_N\www w\Vert^2
\le 2(C\mu^2+1)\Vert \www w\Vert^2.
$$
Moreover, by (1.4), we see $(\www A\www w\cdot \www w) 
\ge \mu\delta\Vert \www w\Vert^2$.
Hence, (2.8) yields
$$
-2s \sum_{k=0}^d (H^k\ppp_k\www w \, \cdot\, \www w)
+ s^2(2\delta\mu-C\mu^2)\Vert \www w\Vert^2
\le \left( 1 + \frac{1}{\ep^2}\right)\Vert Fe^{s\www\va}\Vert^2 
+ 2s^2\ep^2(C\mu^2+1)\Vert \www w\Vert^2.                \eqno{(2.9)}
$$
Since $H^k$, $k=0,1,..., d$ are symmetric, by (2.1) we have
$$
-\sumd (H^k\ppp_k\www w\, \cdot\, \www w)
= -\frac{1}{2} \sumd \int_Q \ppp_k(H^k\www w \, \cdot\, \www w) dxdt
+ \frac{1}{2}\sumd \int_Q ((\ppp_kH^k)\www w\, \cdot\, \www w) dxdt
                                                    \eqno{(2.10)}
$$
\begin{align*}
= &-\frac{1}{2}\int_{\ppp\OOO \times (0,T)}
\sumd ((H^k\nu_k)\www w\, \cdot\, \www w) dSdt 
+ \frac{1}{2}\sumd \int_Q ((\ppp_kH^k)\www w\, \cdot\, \www w) dxdt\\
\ge& -\frac{1}{2}\int^T_0 \int_{\ppp\OOO\setminus \ppp\OOO_-(t)}
\sumd \vert H^k\nu_k\vert \vert \www w\vert^2 dSdt\\
+ &\frac{1}{2}\int^T_0\int_{\ppp\OOO_-(t)} \left\vert
\sumd ((H^k\ppp_k)\www w \cdot \www w)\right\vert dSdt 
- C\Vert \www w\Vert^2.
\end{align*}
Moreover,
$$
-(H^0\ppp_t\www w \,\cdot\, \www w)
= -\frac{1}{2}\int_Q \ppp_t(H^0\www w\cdot \www w) dxdt 
+ \frac{1}{2}\int_Q ((\ppp_tH^0)\www w\cdot \www w) dxdt
                                                 \eqno{(2.11)}
$$
\begin{align*}
=& \frac{1}{2}\int_{\OOO} (H^0(x,0)\www w(x,0) \, \cdot\, \www w(x,0)) dx
- \frac{1}{2}\int_{\OOO} (H^0(x,T)\www w(x,T) \, \cdot\, \www w(x,T)) dx\\
+ &\frac{1}{2}\int_Q ((\ppp_tH^0)\www w\cdot \www w) dxdt\\
\ge & \frac{1}{2}\int_{\OOO} (H^0(x,0)\www w(x,0) \, \cdot\, \www w(x,0)) dx
- \frac{1}{2}\int_{\OOO} (H^0(x,T)\www w(x,T) \, \cdot\, \www w(x,T)) dx
- C\Vert \www w\Vert^2.
\end{align*}
Applying (2.10) and (2.11) in (2.9), we obtain
\begin{align*}
& s^2[(2\delta\mu - C\mu^2) - 2\ep^2(C\mu^2+1)]\Vert \www w\Vert^2
- Cs\Vert \www w\Vert^2\\
+& s\int^T_0 \int_{\ppp\OOO_-(t)} \left\vert 
\sumd ((H^k\nu_k)\www w\, \cdot\, \www w) \right\vert dSdt
+ s\int_{\OOO} (H^0(x,0)\www w(x,0) \, \cdot\, \www w(x,0)) dx\\
\le& \left(1+\frac{1}{\ep^2}\right)\Vert Fe^{s\www\va}\Vert^2
+ Cs\int^T_0\int_{\ppp\OOO\setminus \ppp\OOO_-(t)} \vert \www w\vert^2 dSdt\\
+ & s\int_{\OOO} (H^0(x,T)\www w(x,T) \, \cdot\, \www w(x,T)) dx.
\end{align*}
We choose $\mu>0$ small so that $2\delta\mu - C\mu^2 > 0$.
For this $\mu>0$, we choose $\ep>0$ sufficiently small such that 
$$
(2\delta\mu - C\mu^2) - 2\ep^2(C\mu^2+1) > 0.
$$
Then, for $\mu>0$ and $\ep>0$, we can find a constant $C_0 = C_0(\mu,\ep)>0$ 
such that 
\begin{align*}
& C_0s^2\Vert \www w\Vert^2 - Cs\Vert \www w\Vert^2
+ Cs\int^T_0 \int_{\ppp\OOO_-(t)} \left\vert 
\sumd ((H^k\nu_k)\www w\, \cdot\, \www w) \right\vert dSdt\\
+ & s\int_{\OOO} (H^0(x,0)\www w(x,0) \, \cdot\, \www w(x,0)) dx\\
\le& C_0\Vert Fe^{s\www\va}\Vert^2
+ Cs\int^T_0\int_{\ppp\OOO\setminus \ppp\OOO_-(t)} \vert \www w\vert^2 dSdt
+ s\int_{\OOO} (H^0(x,T)\www w(x,T) \, \cdot\, \www w(x,T)) dx.
\end{align*}
Setting $\www s_0:= \frac{2C}{C_0(\mu,\ep)}$, we obtain
$$
C_0s^2\Vert \www w\Vert^2 - Cs\Vert \www w\Vert^2
\ge \frac{1}{2}C_0s^2\Vert \www w\Vert^2
$$
for all $s> \www s_0$.  Rewriting in terms of $\www w = e^{s\www\va}u
= e^{s\mu\va}u$, we have
\begin{align*}
& \frac{1}{2}C_0s^2\int_Q \vert u\vert^2 e^{2s\mu\va} dxdt
+ s\int^T_0 \int_{\ppp\OOO_-(t)} \left\vert 
\sumd ((H^k\nu_k)u\, \cdot\, u) \right\vert e^{2s\mu\va} dSdt\\
+ & s\int_{\OOO} (H^0(x,0)u(x,0) \, \cdot\, u(x,0)) e^{2s\mu\va(x,0)} dx\\
\le& C_0\Vert F\Vert^2 e^{2s\mu\va} dxdt 
+ Cs\int^T_0\int_{\ppp\OOO\setminus \ppp\OOO_-(t)} \vert u\vert^2 
e^{2s\mu\va} dSdt\\
+ & s\int_{\OOO} (H^0(x,T)u(x,T) \, \cdot\, u(x,T)) e^{2s\mu\va(x,T)} dx
\end{align*}
for all $s > \www s_0$.  Therefore, rewriting $\mu\va$ by $\va$ in the above
inequality, we complete the proof of Theorem 1.
$\blacksquare$
\\

{\bf Remark.}  In Step III, the choice of $B$ is not unique.  For example,
we can similarly argue with 
$$
B:= A - \gamma E_N
$$
with large constant $\gamma > 0$.
The auxiliary estimate in Step III for comleting the proof is 
traditional and conventional for establishing $L^2$-estmates for 
partial differetial equations.  
We can refer for example to the descriptions on 
p.8 in Bellassoued and Yamamoto \cite{BY}, p.37 in Komornik \cite{Ko},
p.10 in Yamamoto \cite{Ya}.
\section{Proof of Theorem 2}

It is known that a relevant Carleman estimate produces an $L^2$-estimate for
initial value (e.g., Kazemi and Klibanov \cite{KK}, Klibanov and Timonov 
\cite{KT}).  Here we apply a method in Huang, Imanuvilov and Yamamoto 
\cite{HIY} which simplifies the argument based on Carleman estimate 
in \cite{KK}, \cite{KT}.
\\

{\bf First Step.}

Let $H^k \in C^1(\ooo{Q})$, $k=0,1,..., d$ be symmetric matrices.
We prove
\\
{\bf Lemma 1 (energy estimate for an initial value problem)}.
\\
{\it
We assume (1.7).  Then there exists a constant $C>0$ such that 
$$
\Vert u(\cdot,t)\Vert^2 
+ \int^T_0 \int_{\ppp\OOO_+(\xi)} 
\left( \sumd (H^k\nu_k)u\, \cdot\, u\right) dSd\xi
$$
$$
\le C\left(\Vert u(\cdot,0)\Vert^2 
+ \int^T_0 \int_{\ppp\OOO \setminus \ppp\OOO_+(\xi)} \vert u\vert^2
dSd\xi\right)
                                               \eqno{(3.1)}
$$
for all $u \in H^1(Q)$ satisfying
$$
H^0\ppp_tu + \sumd H^k\ppp_ku + Pu = 0 \quad \mbox{in $Q$}.
                                                              \eqno{(3.2)}
$$
}

{\bf Proof.}
\\
The proof is by the following standard way.  We take the scalar product 
(3.2) with $u$ and integrate in $x \in \OOO$ to have
\begin{align*}
& (H^0(\cdot,\xi)\ppp_{\xi}u(\cdot,\xi)\, \cdot u(\cdot,\xi))
+ \sumd (H^k(\cdot,\xi)\ppp_ku(\cdot,\xi)\, \cdot u(\cdot,\xi))\\
+& (P(\cdot,\xi)u(\cdot,\xi)\, \cdot u(\cdot,\xi)) = 0, \quad
0\le \xi \le t.
\end{align*}
In the same way as (2.10) and (2.11), we see
\begin{align*}
& \frac{1}{2}\ppp_{\xi} \int_{\OOO}
(H^0(x,\xi)u(x,\xi)\, \cdot u(x,\xi)) dx
+ \frac{1}{2} \int_{\ppp\OOO}
\left( \sumd \nu_kH^k(x,\xi)\, \cdot u(x,\xi)\right) dS\\
-& \frac{1}{2}\int_{\OOO} ((\ppp_tH^0(x,\xi))u(x,\xi)\, \cdot u(x,\xi)) dx
- \frac{1}{2}\int_{\OOO} \left(\left(\sumd \ppp_kH^0(x,\xi)\right)u(x,\xi)\, 
\cdot u(x,\xi)\right) dx\\
+ & \int_{\OOO}(P(x,\xi)u(x,\xi)\, \cdot u(x,\xi))dx = 0, \quad
0\le \xi \le t.
\end{align*}
We integrate in $\xi \in (0,t)$ to have
\begin{align*}
& \int_{\OOO}(H^0(x,t)u(x,t)\, \cdot u(x,\xi)) dx
+ \int^t_0 \int_{\ppp\OOO_+(\xi)}
\left( \left(\sumd \nu_kH^k\right)u(x,\xi)\, \cdot u(x,\xi)\right) dSd\xi\\
=& \int_{\OOO}(H^0(x,0)u(x,0)\, \cdot u(x,0)) dx
- \int^t_0 \int_{\ppp\OOO\setminus \ppp\OOO_+(\xi)}
\left( \left(\sumd \nu_kH^k\right)u(x,\xi)\, \cdot u(x,\xi)\right) dSd\xi\\
+& \int^t_0\int_{\OOO} (\ppp_tH^0(x,\xi)u(x,\xi)\, \cdot\, u(x,\xi)) dxd\xi
\\
+ &\int^t_0 \int_{\OOO} \sumd ((\ppp_kH^k)u(x,\xi)\, \cdot u(x,\xi)) 
dxd\xi
- 2\int^t_0 \int_{\OOO} (Pu(x,\xi)\, \cdot u(x,\xi)) dxd\xi, 
\quad 0\le t \le T.
\end{align*}
We set $E(t) = \int_{\OOO} \vert u(x,t)\vert^2 dx$.
By (1.7), $H^k\in C^1(\ooo{Q})$, and $P \in L^{\infty}(Q)$, we obtain
$$
E(t) + \int^t_0 \int_{\ppp\OOO_+(\xi)}
\left( \left(\sumd \nu_kH^k\right)u(x,\xi)\, \cdot u(x,\xi)\right) dSd\xi
                                                             \eqno{(3.3)}
$$
$$
\le CE(0) + C\int^T_0 \int_{\ppp\OOO \setminus \ppp\OOO_+(\xi)}
\vert u(x,\xi)\vert^2 dSd\xi
+ C\int^t_0 E(\xi) d\xi, \quad 0\le t \le T.
$$
In particular, 
$$
E(t) \le CE(0) 
+ C\int^T_0 \int_{\ppp\OOO \setminus \ppp\OOO_+(\xi)}
\vert u(x,\xi)\vert^2 dSd\xi
+ C\int^t_0 E(\xi) d\xi, \quad 0\le t \le T.
$$
The Gronwall inequality yields
$$
E(t) \le CE(0) 
+ C\int^T_0 \int_{\ppp\OOO \setminus \ppp\OOO_+(\xi)}
\vert u(x,\xi)\vert^2 dSd\xi, \quad
0\le t \le T.
$$
Substituting this into (3.3), we reach
\begin{align*}
& \int^t_0 \int_{\ppp\OOO_+(\xi)}
\left( \left(\sumd \nu_kH^k\right)u(x,\xi)\, \cdot u(x,\xi)\right) dSd\xi\\
\le& CE(0) 
+ C\int^T_0 \int_{\ppp\OOO \setminus \ppp\OOO_+(\xi)}
\vert u\vert^2 dSd\xi, \quad 0\le t \le T.
\end{align*}
The proof of Lemma 1 is complete.

{\bf Second Step.}

We complete the proof of Theorem 2.  By (1.8) we can choose $\beta 
\in \left(0,\,\frac{\delta_0}{M}\right)$ which is close to 
$\frac{\delta_0}{M}$ such that 
$$
T > \frac{1}{\beta}(\max_{x\in \ooo{\OOO}}\eta(x) 
- \min_{x\in \ooo{\OOO}}\eta(x)),
$$
that is,
$$
\min_{x\in \ooo{\OOO}}\va(x,0) - \max_{x\in \ooo{\OOO}}\va(x,T) > \delta_2
> 0,                                          \eqno{(3.4)}
$$
where $\delta_2 > 0$ is some constant.

Moreover, by $\beta < \frac{\delta_0}{M}$, (1.6) and (1.7), we see that 
$$
\sum_{k=0}^d ((\ppp_k\va)H^kv \cdot v) 
= -\beta(H^0v\cdot v)
+ \sumd ((\ppp_k\eta)H^kv \cdot v)
\ge (-\beta M + \delta_0)\vert v\vert^2_{\R^N},
$$
which verifies (1.4) with $\delta := \delta_0 - \beta M$.
Therefore, Theorem 1 is applicable with the weight function 
$\va(x,t) := \eta(x) - \beta t$.

With this $\delta_2$, we apply Theorem 1.  In terms of (1.7), the first term
on the left-hand side of (1.5) is replaced by 
$$
s\delta_1\int_{\OOO} \vert u(x,0)\vert^2 e^{2s\va(x,0)} dx 
\ge s\delta_1e^{2s\min_{x\in \ooo{\OOO}} \va(x,0)} \Vert u(\cdot,0)\Vert^2.
$$
We neglect the second and the third terms on the left-hand side of (1.5)
and further replace the second term on the right-hand side by 
$Cse^{C_1s}\Vert u\Vert^2_{L^2(\ppp\OOO\times (0,T))}$ where
$C_1:= \max_{(x,t)\in \ooo{Q}} \va(x,t)$.

Thus
$$
se^{2s\min_{x\in \ooo{\OOO}} \va(x,0)} \Vert u(\cdot,0)\Vert^2
\le Ce^{C_1s}\Vert u\Vert^2_{L^2(\ppp\OOO\times (0,T))}
+ Cse^{2s\max_{x\in \ooo{\OOO}} \va(x,T)} \Vert u(\cdot,T)\Vert^2
$$
for all $s\ge s_0$.  
Applying (3.4), we obtain
$$
\Vert u(\cdot,0)\Vert^2
\le Ce^{C_1s}\Vert u\Vert^2_{L^2(\ppp\OOO\times (0,T))}
+ Ce^{-2s\delta_2} \Vert u(\cdot,T)\Vert^2           \eqno{(3.5)}
$$
for all $s\ge s_0$. 

Lemma 1 implies
$$
\Vert u(\cdot,T)\Vert^2
\le C(\Vert u(\cdot,0)\Vert^2 
+ \Vert u\Vert^2_{L^2(\ppp\OOO\times (0,T))}).
$$
Substituting this into (3.5), we reach
$$
\Vert u(\cdot,0)\Vert^2
\le Ce^{Cs}\Vert u\Vert^2_{L^2(\ppp\OOO\times (0,T))}
+ Ce^{-2s\delta_2} \Vert u(\cdot,0)\Vert^2
+ Ce^{-2s\delta_2}\Vert u\Vert^2_{L^2(\ppp\OOO\times (0,T))}
$$
for all $s\ge s_0$.   Choosing $s>0$ sufficicently large, we can 
absorb the second term on the right-hand side into the 
left-side hand, and so we complete the proof of Theorem 2.
$\blacksquare$
\vspace{0.2cm}
\section*{Acknowledgment}
The second author was supported by Grant-in-Aid for JSPS Fellows
JP20J11497 of Japan Society for the Promotion of Science. 
The third author was supported by Grant-in-Aid  
for Scientific Research (A) 20H00117 of 
Japan Society for the Promotion of Science, 
The National Natural Science Foundation of China
(no. 11771270, 91730303), and the RUDN University 
Strategic Academic Leadership Program.
The essential part was based on the intensive course by the third author
at Tor Vergata University of Rome and he thanks Professor 
Piermarco Cannarsa for offeirng the opportunity.

\end{document}